\documentclass[aos,MSNbibl,nameyear,dvips]{arximspdf}

\doi{10.1214/13-AOS1124} 
\volume{41}
\issue{4}
\pubyear{2013}
\firstpage{1922}
\lastpage{1943}

\makeatletter
\newcommand{\rrVert}{\Vert}
\newcommand{\rrvert}{\vert}
\newcommand{\llVert}{\Vert}
\newcommand{\llvert}{\vert}

\newtheorem{lemma}{Lemma}
\newtheorem{theorem}{Theorem}

\newcommand{\ebf}{{\mathbf e}}
\newcommand{\cbf}{{\mathbf c}}
\newcommand{\Hbf}{{\mathbf H}}

\makeatother

\begin{document}
\begin{frontmatter}

\title{Asymptotic normality of maximum likelihood and its variational
approximation for stochastic~blockmodels\thanksref{T1}}
\runtitle{Variational approximation for stochastic blockmodels}

\thankstext{T1}{Supported in part by NSF of China under award 11171272.}

\begin{aug}
\author[A]{\fnms{Peter} \snm{Bickel}\ead[label=e1]{bickel@stat.berkeley.edu}},
\author[B]{\fnms{David} \snm{Choi}\corref{}\ead[label=e2]{davidch@andrew.cmu.edu}},
\author[C]{\fnms{Xiangyu} \snm{Chang}\ead[label=e3]{xiangyuchang@gmail.com}}
\and
\author[D]{\fnms{Hai} \snm{Zhang}\ead[label=e4]{zhanghai@nwu.edu.cn}}
\runauthor{Bickel, Choi, Chang and Zhang}
\affiliation{University of California, Berkeley, Carnegie Mellon
University, Xi'an Jiaotong University and Northwest University}
\address[A]{P. Bickel\\
Department of Statistics \\
University of California, Berkeley\\
367 Evans Hall \\
Berkeley, California 94710\\
USA\\
\printead{e1}}
\address[B]{D. Choi\\
Heinz College of Public Policy\\
\quad and Management \\
Carnegie Mellon University \\
Pittsburgh, Pennsylvania 15213 \\
USA\\
\printead{e2}}
\address[C]{X. Chang\\
Department of Applied Mathematics \\
Xi'an Jiaotong University \\
Xi'an 710049 \\
China \\
\printead{e3}}
\address[D]{H. Zhang\\
Department of Mathematics\hspace*{14.88pt} \\
Northwest University \\
Xi'an 710069 \\
China \\
\printead{e4}}
\end{aug}

\received{\smonth{7} \syear{2012}}
\revised{\smonth{4} \syear{2013}}

%
\begin{abstract}
Variational methods for parameter estimation are an active research
area, potentially offering computationally tractable heuristics with
theoretical performance bounds. We build on recent work that applies
such methods to network data, and establish asymptotic normality rates
for parameter estimates of stochastic blockmodel data, by either
maximum likelihood or variational estimation. The result also applies
to various sub-models of the stochastic blockmodel found in the
literature.
\end{abstract}

%

\end{frontmatter}

\section{Introduction}

The analysis of network data is an open statistical problem, with many
potential applications in the social sciences [\citet{lazer2009life}] and
in biology [\citet{proulx2005network}]. In such applications, the models
tend to pose both computational and statistical challenges, in that
neither their fitting method nor their large sample properties are well
understood.

However, some results are becoming known for a model known as the
stochastic blockmodel, which assumes that the network connections are
explainable by a latent discrete class variable associated with each
node. For this model, consistency has been shown for profile likelihood
maximization [\citet{bickel2009nonparametric}], a spectral-clustering
based method [\citet{rohe2011spectral}], and other methods as well
[\citet{bickelmethod,channarond2011classification,choi2010stochastic,coja2008partitioning}], under varying assumptions on the sparsity of the
network and the number of classes. These results suggest that the model
has reasonable statistical properties, and empirical experiments
suggest that efficient approximate methods may suffice to find the
parameter estimates. However, formally there is no satisfactory
inference theory for the behavior of classical procedures such as
maximum likelihood under the model, nor for any procedure which is
computationally not potentially NP under worst-case analysis.

In this note, we establish both consistency and asymptotic normality of
maximum likelihood estimation, and also of a variational approximation
method, considering sparse models and restricted sub-models. To some
extent, we are following a pioneering paper of Celisse et al.
[\citet{celisse2011consistency}], in which the dense model was
considered, and consistency was established, but only for a subset of
the parameters.

\section{Preliminaries}

\subsection{Stochastic blockmodel}

We consider a class of latent variable models considered by various
authors [\citet{karrer2011stochastic,latouche2011overlapping,snijders1997estimation,bickel2009nonparametric}], which we describe as
follows. Let $Z=(Z_1,\ldots,Z_n)$ be latent random variables
corresponding to vertices $1,\ldots,n$, taking values in $[K] \equiv
\{1,\ldots, K\}$. We will assume that $K$ is fixed and does not
increase with $n$. Let $\pi$ be a distribution on $[K]$, and let $H$ be
a symmetric matrix in $[0,1]^{K\times K}$. We define the
\emph{complete graph model} (CGM) for $Z,A$, where $A$ is the $n
\times
n$ symmetric 0--1 adjacency matrix of a graph, by its distribution
%
\begin{equation}
\label{sec2eqCGM}\quad f(Z,A) = \Biggl(\prod_{i=1}^n
\pi(Z_i) \Biggr) \Biggl(\prod_{i=1}^n
\prod_{j=i+1}^n H(Z_i,Z_j)^{A_{ij}}
\bigl(1-H(Z_i,Z_j)\bigr)^{1-A_{ij}} \Biggr),
\end{equation}
where we may interpret $H(Z_i,Z_j)$ as $\mathbb{P}(\mathrm{edge}|Z_i,Z_j)$, and $\pi(a)$ as $\mathbb{P}(Z_i=a)$ for $a=1,\ldots,K$.

The \emph{graph model} (GM) is defined by a distribution $g\dvtx \{0,1\}^{n
\times n} \rightarrow[0,1]$, which satisfies $g(A) = \mathbb
{P}(A;H,\pi
)$ and is given by
\[
g(A) = \sum_{z \in[K]^n} f(z,A).
\]
It is data from GM which we assume we observe.

We will allow $(H,\pi)$ to be parameterized by $\theta$ taking values
in some restricted space $\Theta$, so that parametric submodels of the
blockmodel may be considered. We will consider parameterizations of the
form $\theta= (\rho,\phi)$, in which
\[
H_\theta\equiv\rho S_\phi,\qquad \pi_\theta\equiv
\pi_\phi,\qquad \sum_{a,b=1}^K
\pi_\phi(a)\pi _\phi (b) S_\phi(a,b) = 1,
\]
where $\rho> 0$ is a nonnegative scalar; $\phi$ is a Euclidean
parameter ranging over an open set; $S_\phi$ is a symmetric matrix in
$\mathbb{R}^{K \times K}$; and the map $\phi\mapsto(\pi_\phi,
S_\phi
)$ is assumed to be smooth. Let $\lambda= n\rho$. The interpretation
of these parameters is that $\lambda= \mathbb{E}[\mathrm{degree}]$ and
$\rho= \mathbb{P}(A_{ij}=1)$. 
The utility of this parameterization will be to analyze asymptotic
behavior when $\rho\equiv\rho_n \rightarrow0$ while $\phi$ is kept
fixed, as seems reasonable for sparse network settings.

\subsubsection*{Identifiability of the model} We observe that $f$ is
symmetric under permutation of $Z$ and $\theta$; that is, let $\sigma\dvtx
[K]\rightarrow[K]$ denote a permutation of $[K]$, and let $\Pi$
denote its permutation matrix. For $z \in[K]^n$, let $\sigma(z) =
(\sigma(z_1),\ldots,\sigma(z_n))$, and for $\theta\equiv(\pi,H)$, let
$\sigma(\theta) = (\Pi\pi, \Pi H \Pi^T)$. It then holds for any
permutation $\sigma$ that
\[
f(Z,A; \theta) = f\bigl(\sigma(Z), A; \sigma(\theta)\bigr)
\]
and hence
\[
g(A; \theta) = g\bigl(A; \sigma(\theta)\bigr),
\]
showing that when $Z$ is latent, the stochastic blockmodel is
nonidentifiable. Specifically, $\theta\equiv(\pi, H)$ is equivalent
to $\sigma(\theta) \equiv(\pi\Pi, \Pi H \Pi^T)$ for any permutation
$\sigma$. Let $\mathcal{S}_\theta$ denote this equivalence class, which
corresponds to a relabeling of the latent classes $\{1,\ldots,K\}$. By
an estimate $\hat{\theta}$ under the GM blockmodel, we will mean the
equivalence class $\mathcal{S}_{\hat{\theta}}$. By consistency and
asymptotic normality of $\hat{\theta}$, we will mean that $\mathcal
{S}_{\hat{\theta}}$ contains an element $\theta'$ that converges to the
generative $\theta_0$, or has error $r_n (\theta' - \theta_0)$ that is
asymptotically normal distributed for some rate $r_n \rightarrow0$.

In our analysis, we will assume that the generative $H$ has no
identical rows, as we cannot expect to successfully distinguish classes
which behave identically. If $H$ did contain identical rows, then an
additional source of nonidentifiability would exist. Also, the
generative model would be equivalent to a stochastic blockmodel of
smaller order $K$. We do not treat such cases here.

We note that for some restricted submodels, identifiability can be
restored by imposing a canonical ordering of the latent classes
$1,\ldots,K$. For example, the submodel may restrict $H$ so that
$H(Z_i,Z_j)$ depends only on whether $Z_i=Z_j$ or not; this assumption
could reflect homogeneity of the classes, and is explored in
\citet{rohe2012highest}. This submodel is identifiable under ordering of
$\pi$, and the latent structure might be more gracefully described as a
partition, that is, a variable $X \in\{0,1\}^{n \times n}$ satisfying
$X(i,j) = 1$ iff $Z_i = Z_j$. As a second example, the latent classes
could be ordered by decreasing expected degree. If the submodel
restricts the expected degrees to be unique, the submodel is
identifiable; further discussion can be found in
\citet{celisse2011consistency,bickelmethod}.

\subsubsection*{Degree-corrected blockmodels} An interesting class of
submodels, discussed in
\citet{karrer2011stochastic,zhao2013consistency}, are the
``degree-corrected'' blockmodels with $UV$-many classes obtained by
considering $Z_i = (Z_{i1},Z_{i2})$, for $i=1,\ldots,n$, which take
values $(u,v)$; where $u$ takes values $1,\ldots,U$ with probabilities
$\alpha_1,\ldots,\alpha_U$; and given parameters $\gamma_1,\ldots,\gamma
_V \in[0,1]$, $v$ takes values $\gamma_1,\ldots,\gamma_V$ with
probabilities $\beta_1,\ldots,\beta_V$. We will assume $Z_{i1}$ and
$Z_{i2}$ are independent. Additional parameters needed are a $U \times
U$ symmetric matrix of probabilities $G$. We can now define
\[
\mathbb{P}(Z_{i1} = a, Z_{i2} = \gamma_c,
Z_{j1}=b, Z_{j2}=\gamma_d| A_{ij}=1) =
\alpha_a \alpha_b \beta_c
\beta_d \gamma_c \gamma_d G(a,b).
\]
So although this is a $UV$ blockmodel, it has only $U(U+1)/2+(U-1) +
(2V-1)$ parameters. Its interpretation is that there are $U$ subblocks,
but within each subblock, vertices can hierarchically exhibit further
affinities to vertices both within the same block and other blocks,
thus enabling, for instance, distinction between vertices of high
degree and low degree within each block. This distinction is not
block-dependent, resulting in a reduction of parameters.

Many variants are of course possible; for example, one can choose to
have more parameters by having the $(u,v)$ block probabilities be free,
so that the conditional distribution of $Z_{i2}$ dependent on $Z_{i1}$,
or fewer parameters by treating $\alpha(1),\ldots,\alpha(U)$ as known.

\subsubsection*{More general models} The stochastic blockmodel is a special
case of a more general latent variable model, considered by various authors
[\citet{hoff2002latent,bickel2009nonparametric,bollobs2007phase}]. In
this model, the elements of $Z$ take values in a general space
$\mathcal
{Z}$ rather than $[K]$, $\pi$ is a distribution on $\mathcal{Z}$, and
$H$ is replaced by a symmetric map $h\dvtx \mathcal{Z}\times\mathcal{Z}
\rightarrow[0,1]$. The CGM defines a density for $(Z,A)$, with respect
to an appropriate reference measure, and GM satisfies the identity
%
\begin{equation}
\label{sec2eqidentity} \frac{g(A;\theta)}{g_0(A)} = \mathbb{E}_{\theta_0} \biggl[
\frac
{f(Z,A;\theta)}{f_0(Z,A)}\Big| A \biggr],
\end{equation}
where $f_0$ and $g_0$ denote the distribution under the generative
$\theta_0$. This model is considered in
\citet{bickel2009nonparametric} with $\{Z_i\}_{i=1}^n$ assumed i.i.d.
uniform $(0,1)$. In \citet{handcock2007model}, they are a multivariate
mixture of Gaussian with unknown parameters. If we make no restrictions
on $h$, these models are equivalent.

\subsection{Maximum likelihood and variational estimates}\label
{subsecestimates}

For the complete graph blockmodel, maximum likelihood estimation of $H$
and $\pi$ (or of $\theta$) is basically understood. From (\ref
{sec2eqCGM}) it can be seen that the log likelihood expression
decomposes, so that $\pi$ is estimated from $Z$ independently of $A$,
and $H$ is estimated from $A$ conditional on $Z$. We note that it is
possible for the likelihood to have multiple local optima.\vadjust{\goodbreak}

For the GM blockmodel, the maximum likelihood parameter estimate $\hat
{\theta}^{\mathrm{ML}}$ (i.e., the equivalence class $\mathcal{S}_{\mathrm{ML}} \equiv
\mathcal{S}_{\hat{\theta}^{\mathrm{ML}}}$) is given by
\begin{eqnarray*}
\hat{\theta}^{\mathrm{ML}} &=& \arg\max_{\theta} g(A;\theta)
\\
&=& \arg\max_\theta\sum_{z \in[K]^n}
f(z,A,\theta).
\end{eqnarray*}
Two difficulties in computing $\hat{\theta}^{\mathrm{ML}}$ present themselves:
first, multiple local optima in $g$ may exist even if the CGM
likelihood function $f$ is concave in the appropriate parameterization,
as we shall see for the ordinary unrestricted parameterization. Second,
the maximum likelihood estimate involves a generally intractable
marginalization over the latent variable $Z$.

Variational methods attempt to circumvent the second difficulty (while
accepting the first) by introducing an approximate function $J$ for
which local optimization is computationally easier. For the GM
blockmodel, the estimate $\hat{\theta}^{\mathrm{VAR}}$ (i.e., the equivalence
class $\mathcal{S}_{\hat{\theta}^{\mathrm{VAR}}}$) is given by
\begin{eqnarray*}
\hat{\theta}^{\mathrm{VAR}} &=& \arg\max_\theta
\max_{q \in\mathcal{D}} J(q,\theta;A)
\\
&\stackrel{\triangle} {=}& \arg\max_\theta \max
_{q \in\mathcal
{D}} -D(q\|f_{Z|A;\theta}) + \log g(A;\theta).
\end{eqnarray*}
Here $\mathcal{D}$ is the set of all product distributions over
$\mathcal{Z}^n$, with densities denoted by $\prod_{i=1}^n q_i(\cdot)$.
The term $D(\cdot\|\cdot)$ is the Kullback--Leibler divergence,
and $f_{Z|A;\theta}$ is the conditional density of $Z$ given $A$, that
is, $f_{Z|A;\theta}(Z) = \frac{f(Z,A;\theta)}{g(A;\theta)}$. The
Kullback--Leibler divergence is given by
\[
D(q\| f_{Z|A;\theta}) = \sum_{z \in[K]^n} q(z) \log
\frac
{q(z)}{f_{Z|A;\theta}(z)}.
\]

We note that $J$ simplifies to
\begin{eqnarray*}
J(q,\theta;A) & = & \sum_{i=1}^n \sum
_{a=1}^K q_i(a) \bigl[{-\log
q_i(a)} + \log \pi_\theta(a)\bigr]
\\
&&{}+ \sum_{i=1}^n \sum
_{j=i+1}^n \sum_{a=1}^K
\sum_{b=1}^K q_i(a)q_j(b)
\bigl[ A_{ij} \log H_\theta(a,b)\\
&&\hspace*{138.5pt}{} + (1-A_{ij}) \log
\bigl(1 - H_\theta(a,b)\bigr) \bigr].
\end{eqnarray*}
This formula indicates that, at least for the complete
parameterization, a~local optimum to $J$ can be tractably computed for
moderate $n$ and $K$ using the EM algorithm as in \citet
{daudin2008mixture}. In contrast, optimization of $g$ requires a
summation over $[K]^n$ which is generally intractable. However, note
that we have added $n(K-1)$ new parameters.

Intuitively, we expect the variational estimate to approximate the
maximum likelihood estimate when there exists $q \in\mathcal{D}$ which
is close to $f_{Z|A;\theta}$.

We remark that $\max_q \exp(J(q,\theta;A))$ is upper and lower
bounded by
%
\begin{equation}
\label{equpperlower} f(z,A;\theta) \leq\max_q \exp
\bigl(J(q,\theta;A)\bigr) \leq g(A;\theta)
\end{equation}
for any $z \in[K]^n$. To see this, consider that the lower bound is an
equality if $q = \delta_z$, while the upper bound holds due to
nonnegativity of the Kullback--Leibler divergence. 

\subsubsection*{Other estimation problems} Our focus here is on estimation
of the generative $\theta_0$. In other papers, estimation of the
latent $Z$ is considered to be the primary inferential task
[\citet{rohe2011spectral,choi2010stochastic}]. We feel that both tasks
are of interest. For example, if the data $A$ represents a network
observed in its entirety, estimating $Z$ and quantifying its
uncertainty may give insight into the underlying network structure and
the roles of its actors. On the other hand, if $A$ is understood to be
a representative sample of a larger population, whose overall structure
is of interest, estimates of $\theta$ would be preferable.

\section{Results}

\subsection{Asymptotic normality of maximum likelihood under CGM blockmodel}

We first review the asymptotics of the CGM blockmodel with complete
parameterization.

Parameterize $\theta\equiv(\varpi, \nu)$, where $\varpi\in\mathbb
{R}^K$ and $\nu\in\mathbb{R}^{K \times K}$ are the logit of $\pi$ and~$H$, given by
%
\begin{eqnarray}\label{eqmu}
\varpi(a) &=& \log\frac{\pi(a)}{1 - \sum_{b=1}^{K-1} \pi(b)},\qquad
a=1,\ldots,K-1,
\nonumber\\[-8pt]\\[-8pt]
\nu(a,b) &=& \log\frac{H(a,b)}{1 - H(a,b)},\qquad a,b = 1,\ldots,K,\nonumber
\end{eqnarray}
and let $\mathcal{T}$ denote the canonical parameter space $\{\theta\dvtx
\varpi\in\mathbb{R}^{K-1}, \nu\in\mathbb{R}^{K(K+1)/2}\}$. Let
$(Z,A)$ denote data generated by the model, under the generative
parameter $\theta_0$, and let $f_0$ denote $f$ under $\theta_0$. For
the CGM blockmodel, the log likelihood ratio $\Lambda= \log\frac
{f}{f_0}$ is given by
\begin{eqnarray*}
\Lambda(\theta,Z,A) &=& \sum_{a=1}^{K-1}
\biggl[\bigl(\varpi(a) - \varpi _0(a)\bigr)n_a - n \log
\frac{1 + \sum_{a=1}^{K-1} e^{\varpi(a)}}{1 +
\sum_{a=1}^{K-1} e^{\varpi_0(a)}} \biggr]
\\
&&{}+ \frac{1}{2}\sum_{a=1}^K
\sum_{b=1}^K \biggl[\bigl(\nu (a,b) - \nu
_0(a,b)\bigr)O_{ab} - n_{ab} \log
\frac{1 + e^{\nu(a,b)}}{1 + e^{\nu
_0(a,b)}} \biggr],
\end{eqnarray*}
where
\begin{eqnarray*}
n_a &\equiv& n_a(Z) = \sum_{i=1}^n
1\{Z_i=a\}, \qquad n_{ab} \equiv n_{ab}(Z) =
\sum_{i=1}^n \sum
_{j \neq i}^n 1\{Z_i=a,Z_j=b\},
\\
O_{ab}&\equiv& O_{ab}(A,Z) = \sum
_{i=1}^n \sum_{j \neq i} 1
\{Z_i = a, Z_j=b\} A_{ij}.
\end{eqnarray*}
This is an exponential family in $\theta$. The gradient of $\Lambda$
conditioned on $Z$, evaluated at $\theta' \in\mathcal{T}$, is given by
\begin{eqnarray*}
\frac{\partial\Lambda}{\partial\varpi(a)}\bigl(\theta'\bigr) & = & n_a - n \pi
'(a),\qquad a=1,\ldots,K-1,
\\
\frac{\partial\Lambda}{\partial\nu(a,b)}\bigl(\theta'\bigr) &=& O_{ab} -
n_{ab} H'(a,b),\qquad a,b=1,\ldots,K.
\end{eqnarray*}
Using the parameterization $(\pi, H)$, the maximum likelihood estimates
are given by
\[
\hat{\pi}^{\mathrm{CGM}}(a) = \frac{n_a}{n},\qquad a=1,\ldots,K,
\]
and
\[
\hat{H}^{\mathrm{CGM}}(a,b) = \frac{O_{ab}}{n_{ab}},\qquad a,b=1,\ldots,K.
\]
We note that the paramterizations $(\varpi, \nu)$ and $(\pi,\rho,S)$
are both identifiable under the CGM blockmodel.

\begin{lemma} \label{leCGMnormality2}
Assume the generative parameter $\theta_0 \in\mathcal{T}$ satisfies
$(\log n)^{-1} \lambda_0 \rightarrow\infty$, with $\pi_0$ and $S_0$
constant in $n$. It holds that 
%
\begin{eqnarray*}
\sqrt{n} \bigl(\hat{\varpi}^{\mathrm{CGM}} - \varpi_0\bigr) &
\rightarrow& N(0,\Sigma _1),
\\
\sqrt{n\lambda_0}\bigl(\hat{\nu}^{\mathrm{CGM}}-\nu_0
\bigr) &\rightarrow& N(0,\Sigma_2),
\end{eqnarray*}
where $\Sigma_1$ and $\Sigma_2$ are functions of $\theta_0$.
\end{lemma}
\begin{pf}
The log likelihood ratio $\Lambda$ can be decomposed into two terms
which involve $\varpi$ and $\nu$ separately. Asymptotic normality of
$\hat{\pi}^{\mathrm{CGM}}$ and $\hat{\varpi}^{\mathrm{CGM}}$ follows from standard
exponential family theory. It can be seen that
\begin{eqnarray*}
&&
\sqrt{n\lambda_0} \biggl(\frac{\hat{H}^{\mathrm{CGM}}(a,b)}{\rho_0} - \frac
{H_0(a,b)}{\rho_0}
\biggr) \\
&&\qquad= \sqrt{n^2/n_{ab}} \sqrt{n_{ab}\rho
_0} \biggl(\frac{\hat{H}^{\mathrm{CGM}}(a,b)}{\rho_0} - S_0 \biggr)
\\
&&\qquad= \bigl(\pi_0(a)\pi_0(b) +
o_P(1)\bigr)^{-1/2} \sqrt{n_{ab}\rho
_0} \biggl(\frac{\hat{H}^{\mathrm{CGM}}(a,b)}{\rho_0} - S_0 \biggr),
\end{eqnarray*}
which is asymptotically normal by a Lindeberg central limit theorem.
Since $(\hat{\nu}^{\mathrm{CGM}} - \nu_0) =  (\log\frac{\hat
{H}^{\mathrm{CGM}}(a,b)}{\rho_0} - \log\frac{H_0(a,b)}{\rho_0} +
o_P(1)
)$, asymptotic normality of $\hat{\nu}^{\mathrm{CGM}}$ follows by the delta method.
\end{pf}

Let $\Hbf(\theta') = D_\theta^2 \Lambda(Z,A;\theta)|_{\theta=
\theta
'}$ denote the conditional hessian of $\Lambda$ evaluated at $\theta'$.
For all $\theta' \in\mathcal{T}$, $\Hbf(\theta')$ is given by
%
\begin{eqnarray}
\label{eqhessian1}
\frac{\partial^2 \Lambda}{\partial\varpi(a) \,\partial\varpi
(a)}\bigl(\theta '\bigr) &=& n \pi'(a)
\bigl(1-\pi'(a)\bigr),\qquad a=1,\ldots,K-1,
\\
\label{eqhessian2}
\frac{\partial^2 \Lambda}{\partial\varpi(a) \,\partial\varpi
(b)}\bigl(\theta '\bigr) &=& n \pi'(a)
\pi'(b),\qquad a,b=1,\ldots,K-1,
\\
\label{eqhessian3}
\frac{\partial^2 \Lambda}{\partial\nu(a,b) \,\partial\nu
(a,b)}\bigl(\theta'\bigr) &=& n_{ab}
H'(a,b) \bigl(1-H'(a,b)\bigr),\qquad a,b=1,\ldots,K,\hspace*{-20pt}
\end{eqnarray}
with all other terms equal to zero.

\begin{lemma}[(Local asymptotic normality)] \label{leCGMnormality}
For the CGM blockmodel with parameter values $(\varpi_0$, $\nu_0)
\equiv(\pi_0, \rho_n, S_0) \in\mathcal{T}$, it holds uniformly for
any $s,t$ in a compact set that
%
\begin{eqnarray}
\label{eqCGM} \Lambda \biggl( \varpi_0 + \frac{s}{\sqrt{n}},
\nu_0 +\frac
{t}{\sqrt{n^2
\rho_n}} \biggr) &=& s^T Y_1
+ t^T Y_2 - \frac{1}{2} s^T
\Sigma_1 s \nonumber\\[-8pt]\\[-8pt]
&&{}- \frac{1}{2} t^T \Sigma_2 t
+ o_P(1),\nonumber
\end{eqnarray}
where $\Sigma_1$ and $\Sigma_2$ are functions of $\varpi_0$ and $\nu
_0$, and $Y_1, Y_2$ are asymptotically normal distributed with zero
mean and covariances $\Sigma_1$ and $\Sigma_2$, respectively.
\end{lemma}
\begin{pf} By Taylor expansion,
\begin{eqnarray*}
\Lambda \biggl(\varpi_0 + \frac{s}{\sqrt{n}}, \nu_0 +
\frac{t}{\sqrt{n^2 \rho_n}} \biggr) &=& \Lambda(\varpi_0,\nu_0) +
\frac
{1}{\sqrt{n}} s^T \nabla\Lambda_\varpi(
\theta_0)\\
&&{} + \frac{1}{\sqrt {n\lambda_0}} t^T \nabla
\Lambda_\nu(\theta_0)
 + \frac{1}{n} s^T \Hbf_\varpi(
\theta_0) s\\
&&{} + \frac{1}{n^2 \rho
_n} t^T \Hbf_{\nu}(
\theta_0) t + o_P(1),
\end{eqnarray*}
where $\nabla\Lambda_\varpi(\theta_0)$ and $\nabla\Lambda_\nu
(\theta
_0)$ denote the respective components of the gradient of $\Lambda$
evaluated at $\theta_0$, and $ \Hbf_{\varpi}$ and $ \Hbf_{\nu}$ are
given by (\ref{eqhessian1})--(\ref{eqhessian3}) which describe
$\Hbf
(\theta_0)$. By inspection, $\Lambda(\varpi_0,\nu_0)=0$; $ \Hbf
_{\varpi
}/n$ and $ \Hbf_{\nu}/n^2\rho_n$ converge in probability to constant
matrices; and the random vectors $n^{-1/2} \nabla\Lambda_\varpi$ and
$(n\lambda_0)^{-1/2} \nabla\Lambda_\nu$ converge in distribution by
central limit theorem. This establishes (\ref{eqCGM}).
\end{pf}

For submodels where $\theta\mapsto(\varpi,\nu)$ covers a restricted
subset $\Theta\subset\mathcal{T}$, we generally have $\theta_0 -
\hat
{\theta}^{\mathrm{CGM}} = O_P(n^{-1/2})$. However, if $\theta$ is separable into
$(\theta_\pi,\theta_S)$ such that $\pi= \pi_{\theta_\pi}$ and $S =
S_{\theta_S}$, and $\theta_\pi$ and $\theta_S$ are allowed to vary
freely, then $\theta_S$ has error that is asymptotically normal with
the faster rate $\sqrt{n\lambda}$. Independence of the errors in
$\theta
_S$ and $\theta_\pi$ is then also valid as well. 

\subsection{Asymptotic normality of maximum likelihood under GM blockmodel}

Our main result is that for graphs with poly-log expected degree, the
likelihood ratios of the CGM and GM blockmodels are essentially
equivalent with probability tending to 1, so that inference under the
models is essentially equivalent up to the identifiability restrictions
of the GM blockmodel.
%
\begin{theorem} \label{thconsistency}
Let $(Z,A)$ be generated from a blockmodel with $\theta_0 \in\mathcal
{T}$, such that $S_0$ has no identical columns, and $\rho_0 = \rho_n$
satisfies $n \rho_n / \log n \rightarrow\infty$. Then for all
$\theta
\in\mathcal{T}$,
%
\begin{equation}
\label{eqth1} \frac{g}{g_0}(A,\theta) = \max_{\theta' \in\mathcal{S}_\theta}
\frac
{f}{f_0}\bigl(Z,A,\theta'\bigr) \bigl(1 +
\varepsilon_n\bigl(K,\theta'\bigr)\bigr) + \varepsilon
_n\bigl(K,\theta'\bigr),
\end{equation}
where $\sup_{\theta\in\mathcal{T}} \varepsilon_n(K,\theta) = o_P(1)$.
\end{theorem}
Theorem \ref{thconsistency} is proven in the \hyperref[app]{Appendix}, and can be
viewed as the sum of two parts.
\begin{longlist}[(2)]
\item[(1)] In neighborhoods around $(\varpi_0,\nu_0)$, of order
$(n^{-1/2},(n\lambda)^{-1/2})$, both $f/f_0$ and $g/g_0$ are of order 1
and their difference is $o_P(1)$. We show this using methods similar to
\citet{bickel2009nonparametric}, but it may also be deduced from a
general result in
\citet{le1988preservation}; in their terminology, the profile likelihood
estimate is a distinguished statistic.
\item[(2)] In the exterior of neighborhoods as above, both $f/f_0$ and
$g/g_0$ are both $o_P(1)$ on complements of balls around $\mathcal
{S}_{\theta_0}$ and converge uniformly to~0. Unlike the first, this
part does not seem to follow from
\citet{le1988preservation}.
\end{longlist}

Asymptotic normality of $\hat{\theta}^{\mathrm{ML}}$ follows from Theorem \ref
{thconsistency} and Lemma \ref{leCGMnormality2}, as stated in the
following theorem.
%
\begin{theorem} \label{thnormality}
Assuming the conditions of Theorem \ref{thconsistency} and Lemma \ref
{leCGMnormality}, let $\hat{\varpi}^{\mathrm{ML}}, \hat{\nu}^{\mathrm{ML}}$ and
$\hat
{\varpi}^{\mathrm{CGM}}, \hat{\nu}^{\mathrm{CGM}}$ be the corresponding maximum
likelihood estimates over all $\theta\in\mathcal{T}$. It holds that
$\mathcal{S}_{\mathrm{ML}}$ contains an element $\theta'$ satisfying 
%
\begin{eqnarray*}
\varpi' - \hat{\varpi}^{\mathrm{CGM}} &=& o_P
\bigl(n^{-1/2}\bigr),
\\
\nu' - \hat{\nu}^{\mathrm{CGM}} &=& o_P\bigl((n
\lambda_0)^{-1/2}\bigr).
\end{eqnarray*}
\end{theorem}

\begin{pf}
For each\vspace*{1pt} $\theta' \in\mathcal{S}_{\mathrm{ML} }$ it holds that if either
$|\varpi' - \hat\varpi^{\mathrm{CGM}}| \neq o_P(n^{-1/2})$ or $|\nu' - \hat
\nu
^{\mathrm{CGM}}| \neq o_P((n\lambda_n)^{-1/2})$, then by (\ref{eqCGM}) and
consistency of $\hat{\theta}^{\mathrm{CGM}}$,
\[
\Lambda\bigl(\hat{\theta}^{\mathrm{CGM}}; A,Z\bigr) - \Lambda\bigl(
\theta'; A,Z\bigr) = \Omega_P(1).
\]
Thus, we may prove the lemma by establishing the contrapositive. Since
$\hat{\theta}^{\mathrm{ML}}$ and $\hat{\theta}^{\mathrm{CGM}}$, respectively, maximize
$\frac{g}{g_0}$ and $\frac{f}{f_0}$, it follows by Theorem~\ref
{thconsistency} that for some $\theta' \in\mathcal{S}_{\mathrm{ML}}$,
$\llvert \frac{f}{f_0}(Z,A,\theta^{\mathrm{CGM}}) - \frac{f}{f_0}(Z,A,\theta
')\rrvert  = o_P(1)$, implying that $\Lambda(\hat{\theta}^{\mathrm{CGM}}) -
\Lambda(\theta') = o_P(1)$ for some $\theta' \in\mathcal{S}_{\mathrm{ML}}$.
\end{pf}

A parametrized submodel, such as the degree corrected block model as
discussed earlier,
has likelihood $g(A;\varpi(\theta),\nu(\theta))$. Theorem \ref
{thconsistency} applies, and if the mapping $\theta\mapsto(\varpi,\nu
)$ is smooth, then if estimates for the CGM block submodel exist and
are asymptotically normal, their equivalents in the corresponding GM
model will have equivalent behavior, up to the identifiability issues
that we have discussed. Of course, if CGM block submodel estimates do
not exist or are not consistent, this will be inherited by the GM block
submodel estimates as well.

\subsection{Asymptotic normality of variational estimates under GM blockmodel}

We show that same properties that we have established for maximum
likelihood estimates under the GM blockmodel also hold for the more
computable variational likelihood estimates.

Our proof will use a lemma which follows from the main result of
\citet{bickel2009nonparametric}.
%
\begin{lemma}\label{lepnass}
Let $(A,Z)$ be generated by $\theta_0 \equiv(\rho_n,\pi_0,S_0) \in
\mathcal{T}$, such that $n \rho_n /\log n \rightarrow\infty$ and $S_0$
has no identical columns. It holds that
%
\begin{equation}
\label{eqvarlemma} f(A,Z;\theta_0)/g(A;\theta_0) =
1+o_P(1).
\end{equation}
\end{lemma}
\begin{pf}
By exponential family theory, given a nonidentity permutation $\sigma
$, it holds that $\frac{f}{f_0}(A,Z;\sigma(\theta_0)) = o_P(1)$, and
hence that
$\frac{f}{f_0}(A,\sigma^{-1}(Z);\theta_0) = o_P(1)$ as well. As a result,
%
\begin{equation}
\label{eqvar4} \sum_{Z' \in\mathcal{S}_Z, Z' \neq Z} \frac{f(A,Z';\theta
_0)}{g(A;\theta_0)} \leq
\sum_{Z' \in\mathcal{S}_Z, Z' \neq Z} \frac
{f}{f_0}\bigl(A,Z';
\theta_0\bigr) = o_P(1).
\end{equation}

Given $(A,Z)$ generated under $\theta_0$, let $\hat{z}(A)$ denote the
maximum profile likelihood estimate of $Z$, that is, the set $\arg\max_{z} \sup_\theta f(A,z;\theta)$. Let $\mathcal{S}_Z$ denote the set of
all labels $Z'$ such that $Z' = \sigma(Z)$ for some permutation
$\sigma\dvtx[K] \rightarrow[K]$. Theorem \ref{thconsistency} from
\citet{bickel2009nonparametric} states that under the conditions of
this lemma,
\[
\lim\frac{\log\mathbb{P}_0 (\mathcal{S}_Z \neq\hat
{z}(A))}{\lambda
_n} \leq-s_Q(\pi_0,S_0)
< 0.
\]
This implies that $\mathbb{P}(Z \notin\mathcal{S}_{\hat{z}(A)})=o(1)$.
By Markov's inequality, this implies that $\mathbb{P}(Z \notin
\mathcal
{S}_{\hat{z}(A)}|A) = o_P(1)$, which can be rewritten as
%
\begin{equation}
\label{eqvar3} \sum_{Z' \notin\mathcal{S}_Z} f\bigl(A,Z';
\theta_0\bigr)/g(A;\theta_0) = o_P(1).
\end{equation}
Combining (\ref{eqvar3}) and (\ref{eqvar4}) establishes (\ref
{eqvarlemma}).
\end{pf}

Our result for the variational estimates is Theorem \ref{thvar}.
%
\begin{theorem} \label{thvar}
Let $J(\theta;A)$ denote $\max_{q \in\mathcal{D}} \exp[J(q,\theta;A)]$.
Under the conditions of Theorem \ref{thconsistency} and Lemma
\ref{leCGMnormality},
%
\begin{equation}
\label{eqvarclaim} \frac{ J(\theta;A) }{g(A;\theta_0)} = \max_{\theta' \in\mathcal
{S}_{\theta}}
\frac{f}{f_0}\bigl(Z,A,\theta'\bigr) \bigl(1 +
\varepsilon_n\bigl(K,\theta '\bigr)\bigr) +
\varepsilon_n\bigl(K,\theta'\bigr),
\end{equation}
where $\sup_{\theta\in\mathcal{T}} \varepsilon_n(K,\theta) = o_P(1)$.
Hence, the conclusions of Theorem \ref{thnormality} also apply to
$(\hat{\pi}^{\mathrm{VAR}}, \hat{S}^{\mathrm{VAR}})$, the variational likelihood estimates.
\end{theorem}

\begin{pf}
Recall (\ref{equpperlower}) which states that for all $z$,
\[
f(z,A;\theta) \leq\max_q \exp\bigl(J(q,\theta;A)\bigr)
\leq g(A;\theta).
\]
Dividing the lower bound by $f(A,Z;\theta_0)$, which equals
$g(A;\theta
_0)(1+o_P(1))$ by Lemma \ref{lepnass}, yields
\[
\max_{z \in\mathcal{S}_Z} \frac{f}{f_0}(z,A;\theta) \leq
\frac
{J(\theta;A)}{g(A;\theta_0)(1+o_P(1))}.
\]
The identity $\max_{z \in\mathcal{S}_Z} \frac{f}{f_0}(z,A;\theta) =
\max_{\theta' \in\mathcal{S}_\theta} \frac{f}{f_0}(Z,A;\theta')$
thus implies
%
\begin{equation}
\label{eqvar1} \max_{\theta' \in\mathcal{S}_\theta} \frac{f}{f_0}\bigl(z,A;
\theta'\bigr) \leq \frac{J(\theta;A)}{g(A;\theta_0)(1+o_P(1))}.
\end{equation}
Dividing the upper bound by $g(A;\theta_0)$, and applying Theorem \ref
{thconsistency} yields
%
\begin{equation}
\label{eqvar2}\qquad \frac{J(\theta;A)}{g(A;\theta_0)} \leq\frac{g}{g_0}(A;\theta) \leq\max
_{\theta'\in\mathcal{S}_\theta} \frac{f}{f_0}\bigl(Z,A;\theta '\bigr)
\bigl(1+\varepsilon _n\bigl(K,\theta'\bigr)\bigr) +
\varepsilon_n\bigl(K,\theta'\bigr).
\end{equation}
Combining (\ref{eqvar1}) and (\ref{eqvar2}) to upper and lower bound
$\frac{J(\theta;A)}{g(A;\theta_0)}$ proves the theorem.
\end{pf}

\section{Some statistical applications}

With these results, we can show that some standard inference is valid
using the likelihood or variational likelihood for blockmodels.

\subsection*{Confidence regions for $\theta$}
We have that $\hat{\theta}^{\mathrm{VAR}}$ under $P_{\theta_0}$ is
asymptotically normal with mean $\theta_0$ and variance--covariance
matrices given by Theorem~\ref{thnormality} and Lemma \ref{leCGMnormality2}. Since $\theta\mapsto\Sigma(\theta)$ is continuous, we
can evidently form tests and confidence regions based on $\sqrt {n}(\hat
{\varpi}^{\mathrm{VAR}} - \varpi_0)^T \hat{\Sigma}_1^{-1/2}$ and $\sqrt {n\hat
{\lambda}}(\hat{\nu}^{\mathrm{VAR}} - \nu_0)^T \hat{\Sigma}_2^{-1/2}$, where
$\hat{\Sigma}_1$ and $\hat{\Sigma}_2$ are plug-in estimates of
$\Sigma
_1$ and $\Sigma_2$ using\vspace*{1pt} $\hat{\theta}^{\mathrm{VAR}}$, and $\hat{\lambda}$
equals the average degree in the observed data. The same applies to~$\hat{\theta}^{\mathrm{ML}}$.

\subsection*{Wilks statistic for hypothesis testing}

Under the CGM blockmodel with generative parameter $\theta_0$, the
Wilks (or likelihood ratio) statistic is given by
\[
\Lambda\bigl(Z,A;\hat{\theta}^{\mathrm{CGM}}\bigr) \equiv2\log\frac{f}{f_0}
\bigl(Z,A,\hat {\theta}^{\mathrm{CGM}}\bigr) \rightarrow\chi^2_{{K(K+3)}/{2} - 1}.
\]
This statistic can be used to test against a notional value for $\theta_0$.

A consequence of Theorem \ref{thconsistency} is that
\[
\sup_{\theta\in\mathcal{T}} \log\frac{g}{g_0}(A;\theta) = \sup
_{\theta\in\mathcal{T}} \log \biggl(\frac{f}{f_0}(Z,A;\theta ) \biggr)+
o_P(1),
\]
implying that
\[
\Lambda_G\bigl(A;\hat{\theta}^{\mathrm{ML}}\bigr) \equiv2\log
\frac{g}{g_0}\bigl(A;\hat {\theta }^{\mathrm{ML}}\bigr) = \Lambda\bigl(Z,A;
\hat{\theta}^{\mathrm{CGM}}\bigr) + o_P(1),
\]
so that the Wilks statistic for the GM and CGM estimates have the same
asymptotic distribution, enabling tests against a notional $\mathcal
{S}_{\theta_0}$ using the GM likelihood ratio when $Z$ is latent.

A similar result holds for the Wilks statistic of the variational
estimate $\hat{\theta}^{\mathrm{VAR}}$, which may be easier to compute. To see
this, we observe that since $J(\theta;A) \equiv\max_{q\in\mathcal{D}}
\exp[J(q,\theta;A)] \leq g(A;\theta)$, it holds that
\[
\frac{J(\theta,A)}{J(\theta_0,A)} \geq\frac{J(\theta;A)}{g(A;\theta_0)},
\]
so that Theorem \ref{thvar} implies
\[
\frac{J(\theta,A)}{J(\theta_0,A)} \geq\frac{f}{f_0}(Z,A;\theta) \bigl(1+
o_P(1)\bigr) + o_P(1).
\]
To upper bound the same quantity, we observe that
\begin{eqnarray*}
\frac{J(\theta,A)}{J(\theta_0,A)} &\leq& \frac{g(A;\theta
)}{f(Z,A;\theta
_0)}
\\
&=& \frac{g(A;\theta)}{g(A;\theta_0)f(Z,A;\theta_0)g(A;\theta_0)^{-1}}
\\
&=& \frac{g(A;\theta)}{g(A;\theta_0)(1 + o_P(1))},
\end{eqnarray*}
using Lemma \ref{lepnass}. Thus, the arguments used to bound $\Lambda
_G$ also imply
\[
\Lambda_V\bigl(\hat{\theta}^{\mathrm{VAR}}\bigr) \equiv2\log
\frac{J(\hat{\theta
}^{\mathrm{VAR}},A)}{J(\theta_0,A)} =\Lambda\bigl(Z,A;\hat{\theta}^{\mathrm{CGM}}\bigr) +
o_P(1).
\]

\subsection*{Parametric bootstrap}

The parametric bootstrap is also valid for $\hat{\theta}^{\mathrm{VAR}}$. The
algorithm is:
\begin{longlist}[(2)]
\item[(1)] Estimate $\theta$ by $\hat{\theta}^{\mathrm{VAR}}$.
\item[(2)] Generate $B$ graphs of size $n$ according to the blockmodel with
parameter $\hat{\theta}^{\mathrm{VAR}}$, producing $(Z^*_1,A^*_1),\ldots,(Z^*_B,A^*_B)$.
\item[(3)] Fit $A^*_1,\ldots,A^*_B$ by variational likelihood to get $\hat
{\theta}_1^{\mathrm{VAR}*},\ldots,\hat{\theta}_B^{\mathrm{VAR}*}$.
\item[(4)] Compute the variance--covariance matrix of these $B$ vectors and
use it as an estimate of the truth, or similarly, use the empirical
distribution function of the vectors.
\end{longlist}

\begin{theorem}
Under the conditions of Theorem \ref{thnormality}, the parametric bootstrap
distribution of $\sqrt{n}(\hat{\varpi}^{\mathrm{VAR}} - \varpi_0)$ and
$\sqrt {n\lambda}(\hat{\nu}^{\mathrm{VAR}} - \nu_0)$ converges to the Gaussian limits
given by Lemma \ref{leCGMnormality}.
\end{theorem}

\begin{pf}
Without loss of generality we take $B=\infty$, so that we are asking
that when the underlying parameter is $\hat{\theta}^{\mathrm{VAR}}$, the random
law of $\sqrt{n}(\hat{\varpi}^{\mathrm{VAR}*} - \hat{\varpi}^{\mathrm{VAR}})$ and
$\sqrt {n\lambda}(\hat{\nu}^{\mathrm{VAR}*} - \hat{\nu}^{\mathrm{VAR}})$ converges with
$P_{\theta_0}$ probability tending to 1 to the Gaussian limits of
$\sqrt {n}(\hat{\varpi}^{\mathrm{CGM}} - \varpi_0)$ and $\sqrt{n\lambda}(\hat{\nu
}^{\mathrm{CGM}} - \nu_0)$ as generated under~$\theta_0$.

Let $\hat{\varpi}^{\mathrm{CGM}*},\hat{\nu}^{\mathrm{CGM}*}$ have the distribution of the
CG MLE based on the data that we have generated from $P_{\hat{\theta
}^{\mathrm{VAR}}}$. By standard exponential theory such as our Lemma~\ref{leCGMnormality}, we observe that
%
\begin{eqnarray}
\label{eqbootstrap1} \sqrt{n}\bigl(\hat{\varpi}^{\mathrm{CGM}*} - \hat{
\varpi}^{\mathrm{VAR}}\bigr) &\stackrel {P_{\hat
{\theta}^{\mathrm{VAR}}}} {\longrightarrow}& N(0,
\Sigma_1),
\\
\sqrt{n\lambda}\bigl(\hat{\nu}^{\mathrm{CGM}*} - \hat{\nu}^{\mathrm{VAR}}\bigr)
&\stackrel {P_{\hat
{\theta}^{\mathrm{VAR}}}} {\longrightarrow}& N(0,\Sigma_2),
\end{eqnarray}
since the convergence is uniform on contiguous neighborhoods of $\theta
_0$ and the mapping $\theta\rightarrow(\Sigma_1(\theta),\Sigma
_2(\theta))$ is smooth. As Theorem \ref{thvar} implies local
asymptotic normality, a theorem of Le Cam
[\citet{lehmann2005testing}, Corollary 12.3.1] implies that $P_{\hat
{\theta}^{\mathrm{VAR}}} \triangleleft P_{\theta_0}$ with $P_{\theta_0}$
probability tending to $1$, where $\triangleleft$ denotes contiguity.
As a result, Le Cam's first contiguity lemma (stated below) in
conjunction with Theorem \ref{thvar} implies that
\begin{eqnarray*}
\sqrt{n}\bigl(\hat{\varpi}^{\mathrm{CGM}*} - \hat{\varpi}^{\mathrm{VAR}*}\bigr) &
=& o_{P_{\hat
{\theta}^{\mathrm{VAR}}}}(1),
\\
\sqrt{n\lambda} \bigl(\hat{\nu}^{\mathrm{CGM}*} - \hat{\nu}^{\mathrm{VAR}*}\bigr)
&=& o_{P_{\hat
{\theta}^{\mathrm{VAR}}}}(1).
\end{eqnarray*}
Using this result with (\ref{eqbootstrap1}), it follows that
\begin{eqnarray*}
\sqrt{n}\bigl(\hat{\varpi}^{\mathrm{VAR}*} - \hat{\varpi}^{\mathrm{VAR}}\bigr)
&\stackrel {P_{\hat
{\theta}^{\mathrm{VAR}}}} {\longrightarrow}& N(0,\Sigma_1),
\\
\sqrt{n\lambda}\bigl(\hat{\nu}^{\mathrm{VAR}*} - \hat{\nu}^{\mathrm{VAR}}\bigr)
&\stackrel {P_{\hat
{\theta}^{\mathrm{VAR}}}} {\longrightarrow}& N(0,\Sigma_2),
\end{eqnarray*}
establishing the theorem.

For completeness, we state Le Cam's first contiguity lemma as found in
\citet{van2000asymptotic}, Lemma 6.4.

\begin{lemma}
Let $P_n$ and $Q_n$ be sequences of probability measures on measurable
spaces $(\Omega_n, \mathcal{A}_n)$. Then the following statements are
equivalent:
\begin{longlist}[(2)]
\item[(1)]$Q_n \triangleleft P_n$.
\item[(2)] If $dP_n / dQ_n$ converges in distribution under $Q_n$ to $U$
along a subsequence, then $P(U >0) = 1$.
\item[(3)] If $dQ_n / dP_n$ converges in distribution under $P_n$ to $V$
along a subsequence, then $EV=1$.
\item[(4)] For any statistics $T_n\dvtx \Omega_n \mapsto\mathbb{R}^k$: if $T_n
\stackrel{P_n}{\rightarrow} 0$, then $T_n \stackrel{Q_n}{\rightarrow
} 0$.\qed
\end{longlist}
\end{lemma}
\noqed\end{pf}





\section{Conclusions}

In this paper, we have studied stochastic block and extended
blockmodels, such that the average degree tends to $\infty$ at least at
a polylog rate, and the number of blocks $K$ is fixed. We have shown:
\begin{longlist}[(2)]
\item[(1)] Subject to identifiability restrictions, methods of estimation
and parameter testing on maximum likelihood have exactly the same
behavior as the same methods when the block identities are observed,
such that an easily analyzed exponential family model is in force. The
approach uses the methods of \citet{bickel2009nonparametric} slightly
corrected. Unfortunately, computation of the likelihood is as difficult
as the NP-complete computation of modularities, which also yield
parameter estimates that are usable in the same way.
\item[(2)] We also show that the variational likelihood, introduced in this
context by \citet{daudin2008mixture}, has the same properties as the
ordinary likelihood under these conditions; hence, the procedures
discussed above but applied to the variational likelihood behave in the
same way. The variational likelihood can be computed in $\mathcal
{O}(n^3)$ operations, making this a more attractive method.
\end{longlist}
These results easily imply that classical optimality properties of
these procedures, such as achievement of the information bound, hold.

\subsection*{Discussion} A number of major issues still need to be
resolved. Here are some:
\begin{longlist}[(2)]
\item[(1)] Since the log likelihoods studied are highly nonconcave,
selection of starting points for optimization seems critical. The most
promising approaches from both a theoretical and computational point of
view are spectral clustering approaches
[\citet{rohe2011spectral,chaudhuri2012spectral}].
\item[(2)] Blockmodels play the role of histogram approximations for more
complex models of the type considered in
\citet{bickelmethod}, and if observed covariates are added for models
such as those of
\citet{hoff2002latent}. This implies permitting the number of blocks $K$
to increase, which makes perfect classification and classical rates of
parameter estimation unlikely. Issues of model selection and
regularization come to the fore. Some work of this type has been done in
\citet{rohe2011spectral,choi2010stochastic,chatterjee2012matrix}, but
statistical approximation goals are unclear.
\item[(3)] We have indicated that our results for $(\varpi, \nu
)$-parameterized blockmodels also apply to submodels which are
sufficiently smoothly parameterizable. It seems likely that our methods
can also apply to models where there are covariates associated to
vertices or edges.
\end{longlist}

\begin{appendix}\label{app}
\section*{\texorpdfstring{Appendix: Proof of Theorem \lowercase{\protect\ref{thconsistency}}}{Appendix: Proof of Theorem 1}}

We adopt the convention of \citet{bickel2009nonparametric} and let
$\cbf
$ denote $Z$. Recall that $S = H/\rho_n$. Let $\mu_n = n^2 \rho_n$. Let
$L = \sum_{i \neq j} A_{ij}$. For any $\ebf\in[K]^n$, let
\begin{eqnarray*}
n_{ab}(\ebf) &=& \sum_{i=1}^n
\sum_{j \neq i}^n 1\{\ebf_i=a,\ebf
_j=b\}, \qquad n_a(\ebf) = \sum
_{i=1}^n 1\{\ebf_i=a\},
\\
\pi_a(\ebf) &=& n_a(\ebf)/n, \qquad O_{ab}(A,
\ebf) = \sum_{i=1}^n \sum
_{j
\neq i} 1\{\ebf_i = a, \ebf_j=b\}
A_{ij}.
\end{eqnarray*}
Let $|\ebf- \cbf|$ denote $\sum_{i=1}^n 1\{\ebf_i \neq\cbf_i\}$.
Given $\ebf$, define $\bar{\ebf} = \arg\min_{\ebf' \in\mathcal
{S}_\ebf
} |\ebf' - \cbf|$.
Define the confusion matrix $R \in[0,1]^{K \times K}$ by
\[
\bigl[R(\ebf,\cbf)\bigr]\bigl(a,a'\bigr) = \frac{1}{n} \sum
_{i} 1\bigl\{\ebf_i=a, \cbf
_i=a'\bigr\}.
\]
We observe that for fixed $\cbf$, $R$ is constrained to the set
$\mathcal{R} = \{R\dvtx  R\geq0, R^T1 = \pi(\cbf)\}$. Let $RSR^T \equiv
(RSR^T)(\ebf)$ abbreviate $R(\cbf,\ebf) S R^T(\cbf,\ebf)$. Let
$X(\ebf)
= \mu_n^{-1}\*O(A,\ebf) - RSR^T$.

Let $f_n$ denote the full data likelihood of the stochastic blockmodel,
\[
f_n(A,\ebf;\theta) = \prod_{i=1}^n
\pi_\theta(z_i) \prod_{i < j}
H_\theta(\ebf_i,\ebf_j)^{A_{ij}}\bigl(1
-H_\theta(\ebf_i,\ebf_j)\bigr)^{1-A_{ij}}.
\]
Let $Q_n$ denote the likelihood modularity
[\citet{bickel2009nonparametric}], defined as $Q_n(A,\ebf) = \sup_\theta
\log f_n(A,\ebf;\theta)$. We observe that $Q_n$ equals
\begin{eqnarray*}
Q_n(A,\ebf) &=& \sum_{i=1}^n
n_a \log\frac{n_a}{n}
\\
&&{} + \frac{1}{2} \sum_{a=1}^K \sum
_{b=1}^K \biggl[O_{ab} \log
\frac
{O_{ab}}{n_{ab}} + \biggl(n_{ab}-\frac{O_{ab}}{n_{ab}} \biggr) \log
\biggl(1 - \frac{O_{ab}}{n_{ab}} \biggr) \biggr].
\end{eqnarray*}
For $\rho_n \rightarrow0$, it is shown in
\citet{bickel2009nonparametric} that
\[
Q_n(A,\ebf) = \mu_n \biggl(F \biggl(\frac{O(A,\ebf)}{\mu_n},
\pi (\ebf) \biggr) + \frac{L}{\mu_n}\log\rho_n +
o_P(1) \biggr),
\]
where the function $F$ is given by
\begin{eqnarray*}
F(M,t) &=& \sum_{a,b} t_a
t_b \tau \biggl(\frac{M_{ab}}{t_at_b} \biggr),
\\
\tau(x) &=& x \log x - x
\end{eqnarray*}
for $M \in\mathbb{R}^{K \times K}$ and $t$ in the $K$-simplex.

The result of \citet{bickel2009nonparametric} establishes that the
following properties hold for $F_n$ [also see
\citet{zhao2013consistency} for a reworked derivation]:
\begin{longlist}[(2)]
\item[(1)] The function $\ebf\mapsto F((RSR^T)(\ebf),\pi(\ebf))$ is
maximized by any $\ebf\in\mathcal{S}_\cbf$.
\item[(2)] The function $F$ is uniformly continuous if $M$ and $t$ are
restricted to any subset bounded away from 0.
\item[(3)] Let $G(R,S) = F(RSR^T,R^T1)$. Given $(\pi,S) \in\mathcal{T}$, it
holds for all $R \in\{R \geq0, R^T1=\pi\}$ that
\[
\frac{\partial G((1-\varepsilon)\operatorname{diag}(\pi) + \varepsilon
R,S)}{\partial\varepsilon} \bigg|_{\varepsilon=0+} < -C < 0.
\]
\item[(4)] The directional derivatives
\[
\frac{\partial^2 F}{\partial\varepsilon^2}\bigl(M_0 + \varepsilon(M_1-M_0),
t_0 + \varepsilon(t_1 - t_0)
\bigr)\bigg|_{\varepsilon=0^+}
\]
are continuous in $(M_1,t_1)$ for all $(M_0,t_0)$ in a neighborhood
of\break
$(\operatorname{diag}(\pi)S \operatorname{diag}(\pi),\pi)$.
\end{longlist}

We will use an Bernstein inequality result, similar to that shown in
\citet{bickel2009nonparametric}.
%
\begin{lemma}\label{lepnas}
Let $C_S = \max_{ab} S_{ab}$.
%
\begin{equation}
\label{eqpnaseq1} \mathbb{P} \Bigl(\max_{\ebf} \bigl\|X(\ebf)
\bigr\|_\infty\geq\varepsilon \Bigr) \leq2K^{n+2}\exp \biggl(-
\frac{1}{4}\varepsilon^2 \mu_n \biggr)\vadjust{\goodbreak}
\end{equation}
for $\varepsilon\leq3$, and
%
\begin{eqnarray}
\label{eqpnaseq2}
&&
\mathbb{P} \Bigl(\max_{\ebf:|\ebf-\cbf|\leq m} \bigl\|X(\ebf) -
X(\cbf )\bigr\| _\infty\geq\varepsilon \Bigr)\nonumber\\[-8pt]\\[-8pt]
&&\qquad \leq2\pmatrix{n \cr m}
K^{m+2} \exp \biggl(-\frac
{n}{m(8C_S+2)}\varepsilon^2
\mu_n \biggr)\nonumber
\end{eqnarray}
for $\varepsilon\leq\frac{3m}{n}$.
\end{lemma}
\begin{pf}
$\mu_n X_{ab}$ is a sum of independent zero mean random variables
bounded by $1$. Thus by a Bernstein inequality,
\[
\mathbb{P} \bigl( \bigl|\mu_n X_{ab}(\ebf)\bigr| \geq\varepsilon
\mu_n \bigr) \leq \exp \biggl(\frac{-\varepsilon^2 \mu_n^2}{2(\operatorname{Var}(\mu_n
X_{ab})+ \varepsilon\mu_n/3)} \biggr).
\]
We may bound $\operatorname{Var}(\mu_n X_{ab}) \leq\mu_n$ and
$\varepsilon\leq3$ to yield for fixed $a,b,\ebf$ that
\[
\mathbb{P} \bigl( \bigl|X_{ab}(\ebf)\bigr| \geq\varepsilon\mu_n \bigr)
\leq \exp \biggl(\frac{-\varepsilon^2 \mu_n}{4} \biggr).
\]
A union bound establishes (\ref{eqpnaseq1}).

Similarly, $\mu_n (X_{ab}(\cbf) - X_{ab}(\ebf))$ is a sum of
independent zero mean random variables bounded by $1$. Thus,
\begin{eqnarray*}
&&
\mathbb{P} \bigl( \bigl|\mu_n \bigl(X_{ab}(\ebf)-
X_{ab}(\cbf)\bigr)\bigr| \geq \varepsilon\mu _n \bigr) \\
&&\qquad\leq \exp
\biggl(\frac{-\varepsilon^2 \mu_n^2}{2(\operatorname{Var}(\mu
_n(X_{ab}(\cbf)-X_{ab}(\ebf)))+ \varepsilon\mu_n/3)} \biggr).
\end{eqnarray*}
We may bound $\operatorname{Var}(\mu_n(X_{ab}(\ebf)-X_{ab}(\cbf)))
\leq
4mn C_S \rho_n = 4C_S \mu_n m /n$ and $\varepsilon\leq3m/n$ to yield for
fixed $a,b,\ebf$ that
\[
\mathbb{P} \bigl( \bigl|\mu_n \bigl(X_{ab}(\ebf)-
X_{ab}(\cbf)\bigr)\bigr| \geq \varepsilon\mu _n \bigr) \leq\exp
\biggl(\frac{-\varepsilon^2 \mu_n}{(8C_S +
2)m/n} \biggr).
\]
A union bound establishes (\ref{eqpnaseq2}), where we use that $|\{
\ebf\dvtx  |\ebf-\cbf|\leq m\}| \leq{n \choose m} K^m$ for fixed $\cbf$.
\end{pf}

\begin{pf*}{Proof of Theorem \ref{thconsistency}}
The proof can be separated into four parts.

\subsection*{Part 1: $\ebf$ for which $F$ is small} Here we show, for
some $\delta_n \rightarrow0$, that $F(O(\ebf)/\mu,\pi(\ebf))$ is
suboptimal by at least $\delta_n/2$ for all $\ebf$ in a set
$E_{\delta
_n}$. This will imply that $\sum_{\ebf\in E_{\delta}} \sup_\theta
f(A,\ebf;\theta) = o_P(1) \sup_\theta f(A,\cbf;\theta)$.

By (\ref{eqpnaseq1}), $\mu_n^{-1}O(\ebf) \stackrel{P}{\rightarrow}
RSR^T(\ebf)$ uniformly over $\ebf$; hence, by continuity of $F$ there
exists $\delta_n \rightarrow0$ such that
\[
\mathbb{P} \biggl( \max_{\ebf} \biggl\llvert F \biggl(
\frac{O(\ebf)}{\mu
_n}, \pi (\ebf) \biggr) - F\bigl(RSR^T(\ebf),\pi(\ebf)
\bigr)\biggr\rrvert \geq\delta _n/2 \biggr) = o(1).
\]
As a result, given the sets
\[
E_{\delta_n} = \bigl\{\ebf\dvtx  \bigl\llvert F\bigl(\bigl(RSR^T
\bigr) (\ebf),\pi(\ebf)\bigr) - F\bigl(\bigl(RSR^T\bigr) (\cbf),\pi(
\cbf)\bigr)\bigr\rrvert \geq\delta_n \bigr\},
\]
it holds\vspace*{1pt} for all $\ebf\in E_{\delta_n}$ that $F (\frac{O(\ebf
)}{\mu
_n},\pi(\ebf) ) \leq F(RSR^T(\cbf),\pi(\cbf)) - \delta_n/2 +
o_P(\delta_n)$. We may choose $\delta_n$ to additionally satisfy
%
\begin{eqnarray}\label{eqpart1}
\sum_{\ebf\in E_\delta} e^{\mu_n F((RSR^T)(\ebf),\pi
(\ebf))} &\leq& \sum
_{\ebf\in E_\delta} e^{\mu_n  (F((RSR^T)(\cbf),\pi
(\cbf
)) + o_P(\delta_n) - \delta_n/2 )}
\nonumber\\
&\leq& e^{\mu_n F((RSR^T)(\cbf),\pi(\cbf))} e^{-\mu_n
(1+o_P(1)) \delta_n/2}K^n
\\
&=& e^{\mu_n F((RSR^T)(c),\pi_c)} o_P(1),\nonumber
\end{eqnarray}
where we require $\delta_n\rightarrow0$ slowly enough that $\mu_n
\delta_n \gg n$.

\subsection*{Part 2: A concentration inequality} We wish to show for
$\ebf\notin E_{\delta_n}$ a result similar to part 1. However, as some
$\ebf$ will be very close to $\cbf$, we must bound the suboptimality of
$F(O(\ebf)/\mu,\pi(\ebf))$ more carefully.

By (\ref{eqpnaseq2}), it holds that
\[
\mathbb{P} \biggl(\max_{\ebf:|\ebf-\cbf|=m} \bigl\|X(\ebf)-X(\cbf)\bigr\|
_\infty \geq\varepsilon\frac{m}{n} \biggr) \leq2n^mK^{m+2}
\exp \biggl(-\frac
{m}{n(8C_S+2)}\varepsilon^2\mu_n \biggr).
\]
It follows that we may choose $\varepsilon\rightarrow0$ such that
\begin{eqnarray*}
\mathbb{P} \biggl(\max_{\ebf\notin\mathcal{S}_\cbf} \frac{\|
X(\bar{\ebf
})-X(\cbf)\|_\infty}{|\bar{\ebf}-\cbf|/n} \geq\varepsilon
\biggr) &\leq& \sum_{m=1}^n \mathbb{P}
\biggl(\max_{\ebf:|\ebf-\cbf|=m, \ebf
=\bar{\ebf}} \frac{\|X(\ebf)-X(\cbf)\|_\infty}{m/n} \geq\varepsilon \biggr)
\\
&\leq& \sum_{m=1}^n 2 K^K
n^mK^{m+2}\exp \biggl(-\frac
{m}{n(8C_S+2)}
\varepsilon^2\mu_n \biggr)
\\
&\leq& \sum_{m=1}^n 2K^{K+2}e^{m (\log n + \log K -
{\varepsilon
^2 \mu_n }/({n(8C_2+2)}) )}
\\
&=& o(1), 
\end{eqnarray*}
where the final equality holds because $\mu_n/n \gg\log n$, so that we
may choose $\varepsilon\rightarrow0$ such that $\varepsilon^2 \mu_n/n \gg
\log n$. It follows that
%
\begin{equation}
\label{eqpart2} \max_{\ebf\notin\mathcal{S}_\cbf} \frac{\|X(\bar{\ebf}) -
X(\cbf)\|_\infty}{|\bar{\ebf}-\cbf|/n} = o_P(1).
\end{equation}

\subsection*{Part 3: $\ebf$ when $F$ is large} Here we bound the
suboptimality of $F(O(\ebf)/\break\mu,\pi(\ebf))$ in similar fashion to
part 1.

Recall $F((RSR^T)(\ebf),\pi(\ebf)) = G(R(\ebf),S)$ with $R(\ebf)
\in\{
R \geq0,R^T1 = \pi(\cbf)\}$. Let $h(\ebf)$ abbreviate $RSR^T(\ebf) -
RSR^T(\cbf)$. Property 3 implies that for all $\ebf$,
\[
\frac{\partial}{\partial\varepsilon} F \bigl(RSR^T(\cbf) +\varepsilon h(\ebf ),
\bigl(RSR^T(\cbf) +\varepsilon h(\ebf) \bigr)^T1 \bigr)
\bigg|_{\varepsilon
=0^+} < -\Omega_P(1),
\]
where $a_n = \Omega_P(b_n)$ denotes that $a_n$ is bounded below (in
probability) by $b_n$ times a constant factor. As $\delta_n
\rightarrow
0$, this implies for all $\ebf\notin E_{\delta_n}$,
\[
F \bigl(\bigl(RSR^T\bigr) (\cbf),\pi(\cbf) \bigr) - F \bigl(
\bigl(RSR^T\bigr) (\ebf),\pi (\ebf ) \bigr) \geq\frac{1}{n}
\Omega\bigl(|\bar{\ebf} - \cbf|\bigr).
\]
As $(O(\cbf)/\mu_n, \pi(\cbf))$ converges in probability to
$(RSR^T(\cbf
),\pi(\cbf))$, properties 3 and~4 together imply for all $\ebf$,
\[
\frac{\partial}{\partial\varepsilon} F \biggl(\frac{O(\cbf)}{\mu_n} +\varepsilon h(\ebf), \biggl(
\frac{O(\cbf)}{\mu_n} +\varepsilon h(\ebf) \biggr)^T1 \biggr)
\bigg|_{\varepsilon=0^+} < -\Omega_P(1)
\]
and thus for $\ebf\notin E_{\delta_n}$,
\[
F \biggl(\frac{O(\cbf)}{\mu_n}, \pi(\cbf) \biggr) - F \biggl(\frac
{O(\cbf
)}{\mu_n} +
h(\bar{\ebf}), \pi(\bar{\ebf}) \biggr) \geq\frac
{1}{n}\Omega
_P\bigl(|\bar{\ebf}-\cbf|\bigr)
\]
and hence also that
%
\begin{equation}
\label{eqlastrevision}\qquad F \biggl(\frac{O(\cbf)}{\mu_n}, \pi(\cbf) \biggr) - F
\biggl(\frac
{O(\cbf
)}{\mu_n} + h(\bar{\ebf}) \bigl(1+o_P(1)\bigr), \pi(
\bar{\ebf}) \biggr) \geq \frac
{1}{n}\Omega_P\bigl(|\bar{\ebf}-
\cbf|\bigr).
\end{equation}
It can be seen that $h(\bar{\ebf}) \equiv RSR^T(\bar{\ebf}) -
RSR^T(\cbf ) = \Omega (\|\bar{\ebf} - \cbf\|/n )$. As a result,
by~(\ref{eqpart2}), for all $\ebf\notin E_{\delta_n}$,
\begin{eqnarray*}
&&
\biggl\llVert \frac{O(\bar{\ebf})}{\mu_n} - \frac{O(\cbf)}{\mu_n} - \bigl(RSR^T(
\bar{\ebf}) - RSR^T(\cbf) \bigr) \biggr\rrVert _\infty\\
&&\qquad=
o_P \bigl(|\bar{\ebf}-\bar{\cbf}|/n \bigr)
\\
&&\qquad= o_P \bigl(RSR^T(\bar{\ebf}) - RSR^T(
\cbf) \bigr)
\end{eqnarray*}
and hence manipulation yields for all $\ebf\notin E_{\delta_n}$,
%
\[
\frac{O(\bar{\ebf})}{\mu_n} - \frac{O(\cbf)}{\mu_n} = \bigl(RSR^T(\bar {\ebf}) -
RSR^T(\cbf) \bigr) \bigl(1+o_P(1)\bigr),
\]
where the $o_P(1)$ term is uniform over $\ebf$. As a result, it follows
from (\ref{eqlastrevision}) that for $\ebf\notin E_{\delta_n}$,
\[
F \biggl(\frac{O(\cbf)}{\mu_n},\pi(\cbf) \biggr) - F \biggl(\frac
{O(\ebf
)}{\mu_n},
\pi(\bar{\ebf}) \biggr) \geq\frac{1}{n}\Omega_P\bigl(|\bar {\ebf }-
\cbf|\bigr),
\]
where the $\Omega_P(|\bar{\ebf}-\cbf|)$ is uniform over $\ebf$. It
follows that
%
\begin{eqnarray}\label{eqpart3}
&&
\sum_{\ebf\notin E_\delta, \ebf\notin\mathcal{S}_\cbf}
e^{\mu_n F ({O(\ebf)}/{\mu_n},\pi(\ebf
) )} \nonumber\\
&&\qquad\leq \sum_{m=1}^n \sum
_{\ebf:|\bar{\ebf}-\cbf|=m}
e^{\mu_n
F ({O(\ebf)}/{\mu_n},\pi(\ebf) )}
\nonumber\\
&&\qquad= \sum_{m=1}^n \sum
_{\ebf:|\bar{\ebf}-\cbf|=m}
e^{\mu_n  [ F ({O(\cbf)}/{\mu_n},\pi(\cbf) )
+ F
({O(\ebf)}/{\mu_n},\pi(\ebf) )-F ({O(\cbf
)}/{\mu_n},\pi
(\cbf) ) ]}
\nonumber\\
&&\qquad\leq \sum_{m=1}^n \sum
_{\ebf:|\bar{\ebf}-\cbf|=m}
e^{\mu_n F ({O(\cbf)}/{\mu_n},\pi(\cbf) )} e^{-
{\mu_n}\Omega_P(m)/{n}}
\\
&&\qquad\leq \sum_{m=1}^n
e^{\mu_n F ({O(\cbf)}/{\mu
_n},\pi
(\cbf) )} K^K n^m K^m e^{-{\mu_n}\Omega_P(m)/{n}}
\nonumber\\
&&\qquad\leq \sum_{m=1}^n
e^{\mu_n F ({O(\cbf)}/{\mu
_n},\pi
(\cbf) )} K^K e^{m (\log n + \log K - \Omega_P(\mu
_n/n)
)}
\nonumber\\
&&\qquad = e^{\mu_n F ({O(\cbf)}/{\mu_n},\pi(\cbf) )} o_P(1).\nonumber
\end{eqnarray}
%



\subsection*{Part 4: Putting the parts together}
Combining (\ref{eqpart3}) and (\ref{eqpart1}) yields that
%
\begin{equation}
\label{eqpart4-1} \sum_{\ebf\notin\mathcal{S}_\cbf} e^{\mu_n F ({O(\ebf
)}/{\mu
},\pi(\ebf) )} \leq
e^{\mu_n F((RSR^T)(\cbf),\pi(\cbf))} o_P(1).
\end{equation}
Since $\frac{f}{f_0}(A,\cbf;\theta)$ is unimodal in $\theta$, it holds
that if $\frac{f}{f_0}(A,\cbf;\theta) \neq o_P(1)$, then $\theta
\rightarrow\theta_0$, and hence by Lemma \ref{leCGMnormality2},
$\frac{f}{f_0}(A,\cbf;\sigma(\theta)) = o_P(1)$ for any nonidentity
permutation $\sigma$. It follows that
%
\begin{eqnarray}\label{eqpart4-2}
\sum_{\cbf' \in\mathcal{S}_\cbf} f\bigl(A,\cbf';
\theta\bigr) &=& \sum_{\theta' \in\mathcal{S}_\theta} f\bigl(A,\cbf;
\theta'\bigr)
\nonumber\\[-8pt]\\[-8pt]
& = & \max_{\theta' \in\mathcal{S}_\theta} f(A,\cbf;\theta) \bigl(1+o_P(1)
\bigr).\nonumber
\end{eqnarray}
Combining (\ref{eqpart4-1}) and (\ref{eqpart4-2}) yields
\[
\sum_{\ebf\neq\cbf} \sup_\theta f(A,\ebf;
\theta) = \Bigl(\sup_\theta f(A,\cbf;\theta)
\Bigr)o_P(1).
\]
Letting $F_0$ abbreviate $\sup_\theta f(A,\cbf;\theta)$, and using
$g(A;\theta) = \sum_{\ebf} f(A,\ebf;\theta)$,
\begin{eqnarray*}
\frac{g(A;\theta)}{g(A;\theta_0)} &=& \frac{\sum_{\ebf} f(A,\ebf;\theta
)}{\sum_{\ebf} f(A,\ebf;\theta_0)}
\\
&=& \frac{ f(A,\cbf;\theta)}{f(A,\cbf;\theta_0) + \sum_{\ebf\neq
\cbf}
f(A,\ebf;\theta_0)} \\
&&{}+ \frac{ \sum_{\ebf\neq\cbf} f(A,\ebf;\theta
)}{f(A,\cbf;\theta_0) + \sum_{\ebf\neq\cbf} f(A,\ebf;\theta
_0)}
\\
&=& \frac{ f(A,\cbf;\theta)}{f(A,\cbf;\theta_0) + F_0 o_P(1)} + \frac{
F_0 o_P(1) }{f(A,\cbf;\theta_0) + F_0 o_P(1)},
\end{eqnarray*}
where in the last equality we have used the fact that for all $\theta$
\[
0 \leq\sum_{\ebf\neq\cbf} f(A,\ebf;\theta) \leq\sum
_{\ebf\neq
\cbf
} \sup_\theta f(A,\ebf;\theta) =
F_0 o_P(1).
\]
Since $F_0$ equals the likelihood of the MLE under the CGM model, it
holds that $\frac{F_0}{f(A,\cbf;\theta_0)}$ converges in distribution,
and hence $\frac{F_0}{f(A,\cbf;\theta_0)} o_P(1) = o_P(1)$. We may
therefore substitute $F_0 o_P(1) = f(A,\cbf;\theta_0) \frac
{F_0}{f(A,\cbf;\theta_0)} o_P(1) = f(A,\cbf;\theta_0) o_P(1)$ to yield
\begin{eqnarray*}
\frac{g(A;\theta)}{g(A;\theta_0)} &=& \frac{ f(A,\cbf;\theta
)}{f(A,\cbf;\theta_0)(1 + o_P(1))} + \frac{ f(A,\cbf;\theta_0) o_P(1)
}{f(A,\cbf;\theta_0)(1+ o_P(1))}
\\
&=& \frac{f}{f_0}(A,\cbf;\theta) \bigl(1 + o_P(1)\bigr) +
o_P(1),
\end{eqnarray*}
which proves the theorem.
\end{pf*}
\end{appendix}

\section*{Acknowledgements}

We would like to thank the reviewers for their help in fixing an
earlier version of the paper.



\printaddresses

\end{document}